\documentclass{amsart}
\usepackage{amssymb,amsmath}
\usepackage{amsthm}
\usepackage{euscript}

\usepackage{graphics}

\usepackage[greek, english]{babel}
\usepackage{teubner}

\input{epsf}

\newcounter{punct}

\def\punct{\refstepcounter{punct}{
\arabic{punct}.  }}

\newtheorem{theorem}{Theorem}

\newtheorem{observation}[theorem]{Observation}

          \def\sm{\smallskip}

\newcommand{\Aut}{\mathop {\mathrm {Aut}}\nolimits}

\def\0{\mathbf 0}

\def\wh{\widehat}

\renewcommand{\Aut}{\mathop {\mathrm {Aut}}\nolimits}

 \def\C {{\mathbb C }}
  \def\Z{{\mathbb Z}}

\def\Q{{\mathbb Q}}

\def\kappa{\varkappa}
\def\epsilon{\varepsilon}
\def\phi{\varphi}

\def\Aut{\operatorname{Aut}}

\def\cH{\EuScript{H}}

\def\0{\boldsymbol{0}}

\def\kop{\text{\bf \Koppa}}

\begin{document}

\begin{center}
\bf\large
On algebras of conjugacy classes of symmetric groups
\\ with respect
to Young subgroups
\bigskip

\sc\large
Yury A. Neretin%
\footnote{Supported by the grant FWF (The Austrian Scientific Funds),
Project PAT5335224.}
\end{center}

\bigskip

{\small
We consider the subalgebra of the group algebra of a symmetric group consisting of functions that are constant on conjugacy classes with respect to a Young subgroup. We write an expression for structure constants of this algebra in the spirit of Hurwitz numbers.

}

\bigskip

{\bf\punct Introduction.} Centers of the group algebras of symmetric groups $S_N$ are well-known
important and difficult objects, it was a topic of numerous 
investigations, see for example \cite{Hur}, \cite{FH}, \cite{Ju}, \cite{IK}, \cite{Gri}. In particular, structure constants of such algebras
can be expressed as Hurwitz numbers (certain expressions including 
summations over sets of ramified coverings over a compact two-dimensional
surface), see, \cite{Hur}, \cite{LZ}, \cite{NO}.

Many algebras of conjugacy classes and double cosets%
\footnote{Algebras of conjugacy classes are special cases
of algebras of double cosets.}
 for finite
and locally compact groups (as algebras of double cosets of semisimple
groups with respect to maximal compact subgroups, Hecke algebras, affine Hecke algebras) also were topics of numerous investigations. On the 
other hand, R.~S.~Ismagilov and G.~I.~Olshanski observed
that for infinite-dimensional groups spaces of double cosets can admit
natural semigroup structures. This phenomenon is quite usual,
known picture includes numerous double cosets spaces, which seem too general in classical theory
 (for symmetric groups, see \cite{Olsh}, \cite{Ner-tri},  \cite{Ner-UMN}, \cite{Ner-encoding}). The author did several experiments with descent from infinite-dimensional
 level to finite and 
 finite-dimensional groups, \cite{Ner-young}, \cite{Ner-IK}, \cite{Ner-smaller}. This note is a continuation of such experiments.

\sm

{\bf\punct Young subgroups.} For a finite set $X$ denote by $\# X$ the number of its elements.

 Let $\Omega$ be a disjoint union
$$\Omega=\Omega_1\cup \Omega_2\cup \dots\cup\Omega_\psi$$
of finite sets,
$$
\#\Omega_j=m_j,\quad \#\Omega=N=\sum m_j.
$$
We say, that we have a collection of colors $\gimel_1$, \dots, 
$\gimel_\psi$,
and elements of each $\Omega_j$ are painted
in  a color $\gimel_j$. 

We consider the group $S_N$  of all permutations of $\Omega$
and its {\it Young subgroup}
$$
Y=Y(\{m_j\})=S_{m_1}\times\dots\times S_{m_\psi} 
$$
 consisting of  permutations preserving the coloring, i.e., permutations
of $\Omega$ preserving each subset $\Omega_j$.

\sm

\begin{figure}
$$1)\epsfbox{necklace.1}\quad 2) \epsfbox{necklace.2}\quad 
3)\epsfbox{necklace.3}\quad 4)\epsfbox{necklace.4}
\quad 5)\epsfbox{necklace.5}$$
\\
Five necklaces. The necklaces 3) and 4) are isomorphic.
Groups $\Aut(\cdot)$ are respectively cyclic groups $\Z_1$, $\Z_2$, $\Z_3$, $\Z_3$,
$\Z_1$. If we regard this collection as a neck mess, the group
of automorphisms is $\Z_2\times \bigl(S_2\ltimes (\Z_3)^2\bigr)$ 
\caption{}\label{fig:1}
\end{figure}

\begin{figure}
a)\epsfbox{checker.1}\qquad b)\quad \epsfbox{checker.2}
\\
Notation for colors. NearIR: $\epsfbox{checker.3}$,
 midIR: \epsfbox{checker.4}, farUV: \epsfbox{checker.5}. 
\\
a) Plus-triangles and minus-triangles.
\\
b) A piece of a checker triangulated surface.
\caption{}\label{fig:2}
\end{figure}

\begin{figure}
a)\epsfbox{checker.6} b)\qquad \epsfbox{checker.7}
\\
a) A chamber.
\\
b) Permutations $u$ (arrows \epsfbox{checker.8})
and $v$ (arrows \epsfbox{checker.9}).
\caption{}\label{fig:3}
\end{figure}

\begin{figure}
\epsfbox{vertex.1}\qquad \epsfbox{vertex.2}\qquad \epsfbox{vertex.3}
\\
Cycles around nearIR-farUV, midIR-farUV, nearIR-midIR vertices.
\caption{}\label{fig:4}
\end{figure}

{\bf\punct  Necklaces and neck messes.}
A {\it necklace} $p$ of length $\ell$ is a cyclic collection of $\ell$ beads
 of $\psi$
possible colors defined up to a cyclic permutation
(see Fig. \ref{fig:1}). An isomorphism $p\to q$ of necklaces
 is a bijection between beads
of $p$ and $q$ preserving colors and cyclic order of colors.
An automorphism of $p$ is an isomorphism of $p$ to itself. 
 By $\Aut(p)$ we denote the group of automorphisms of
$p$, i.e. the group of rotations, preserving colors.
Clearly, $\Aut(p)$ is a cyclic group  whose order divides the length
 $\ell$ (see Fig. \ref{fig:1}).

A {\it neck mess} $P=\{p_i\}$ is a finite collection of necklaces.
An isomorphism of two neck messes $P=\{p_1,\dots,p_k\}$, 
$Q=\{q_1,\dots,q_k\}$ is a collection  of isomorphisms
$p_i\to q_{\sigma(i)}$, where $\sigma\in S_k$ is a permutation
of the set $\{1,\dots,k\}$. 

Let a neck mess $P=\{p_i\}$ consists of $\alpha_1$ copies of a necklace
$r_1$, $\alpha_2$ copies of a necklace $r_2$, etc. Let $r_l$ be pairwise
nonisomorphic. Then the group of automorphisms of $P$ is isomorphic
to
$$
\Bigl(S_{\alpha_1}\ltimes \bigl(\Aut(r_1)\bigr)^{\alpha_1}\Bigr)
\times
\Bigl(S_{\alpha_2}\ltimes \bigl(\Aut(r_1)\bigr)^{\alpha_2}\Bigr)
\times \dots,
$$ 
where the symbol $\ltimes$ denotes a semi-direct product,
and $S_{\alpha_l}$ acts on $\bigl(\Aut(r_1)\bigr)^{\alpha_2}$
by permutations of factors (see Fig. \ref{fig:1}).
 
\sm

{\bf\punct Conjugacy classes of $S_N$ with respect to the Young subgroups.}
Consider conjugacy classes of  $S_N$ with respect to the subgroup
$Y=Y(\{m_j\})$, i.e., classes of equivalence on $S_N$ defined by
$$
g\sim hgh^{-1},\quad\text{where $h\in Y$.}
$$
Denote by $S_N/\!/ Y$ the quotient space.
Denote by $\C[S_N/\!/ Y]$ the subalgebra of the group algebra
$\C[S_N]$ consisting of functions, which are constant on such conjugacy
classes.

\begin{observation}
\label{obs:conjugacy}
 Conjugacy classes $S_N/\!/ Y$ are in one-to-one correspondence 
with neck messes having $m_j$ beads of each color $\gimel_j$; neck messes
are defined up to isomorphisms.
\end{observation}

The construction of the correspondence is the following. We take a representative of a conjugacy class and decompose it as a product of disjoint cycles. Each element of a cycle has some color $\gimel_i$,
and therefore each cycle determines a necklace.
\hfill $\square$

\sm

{\bf\punct Checker triangulated surfaces.}
A {\it checker triangulated surface} is a colored graph on an oriented
not necessary connected closed surface satisfying the following conditions:

\sm

--- edges of the graph are painted 3 colors, say near-infrared (nearIR), mid-infrared (midIR), far-ultraviolet (farUV)%
\footnote{We  avoid intersections with colors $\gimel_j$, 
which can used for beads.};

\sm

--- the graph divides the surface into triangles of two types,
shown on 
 Fig. \ref{fig:2}.a,
say plus-triangles and minus-triangles;  

\sm

--- neighbors of plus-triangles are-minus triangles and wise versa. 

\sm

See  Fig. \ref{fig:2}.b. 

We also assign mixed colors for vertices of the triangulation
according types of edges entering to this vertex,
nearIR-farUV, midIR-farUV, nearIR-midIR.

We denote by $\kop_N$ the set of  surfaces with $2N$
triangles defined up to a combinatorial equivalence. 

\sm 

{\sc Remark.} This combinatoric structure arise for various reasons in different domains of mathematics, see \cite{ShV}, \cite{LZ}, 
\cite{Ner-tri}, \cite{Ner-virtual}.
The main origin of the interest to such objects is Belyi theorem
and the related Grothendieck program of investigation of the Galois group
of the rational numbers $\Q$ using graphs on surfaces (dessins d'enfant), see for example \cite{ShV}, \cite{Sha},
\cite{LZ}, \cite{GG}.

 \sm

{\bf\punct Chambers and labelings.} Consider a checker triangulated surface
$\cH\in  \kop_N$. We say that a {\it chamber} of $\cH$ is
a quadrilateral obtained as a union of two triangles having a common
farUV edge, see Fig. \ref{fig:3}.





We say that a {\it labeling} of an element  of $\cH\in \kop_N$ is 
an enumeration of chambers by numbers $1$, \dots, $N$.
Denote by $\wh\kop_N$ the set of all labeled 
 surfaces obtained in this way.

\sm

\begin{observation}
 There is a canonical one-to-one correspondence
$$
S_N\times S_N\leftrightarrows \wh\kop_N.
$$
\end{observation}

Namely, for an element of $\wh\kop_N$ we define
$u$, $v\in S_N$
in the following way. For each nearIR  edge
we assume that $u$ 
 sends the label on  plus-side to the label on
minus-side of the edge. Similarly, 
for each  midIR edge
we assume that $v$ 
 sends the label on  minus-side to the label on
plus-side of the edge, see Fig. \ref{fig:3}.b.

Conversely, consider a pair of elements of $S_N$, say 
$u$, $v$.
Consider a collection of $N$ disjoint labeled chambers.

\sm
 
--- if $u$
sends $\alpha\mapsto\beta$, then we glue together 
$\alpha$-th chamber  and $\beta$-th chamber
along nearIR plus-edge of 
$\alpha$-th chamber and nearIR minus-edge of $\beta$-th chamber.

\sm

--- if $v$
sends $\gamma\mapsto\delta$, then we glue together
$\gamma$-th chamber and $\delta$-th chamber
along midIR minus-edge of 
$\gamma$-th chamber and midIR plus-edge of $\delta$-th chamber.
\hfill $\square$

\sm

{\sc Remark.} See a big collection of constructions of such type in
\cite{Ner-encoding}.
\hfill $\boxtimes$

\sm

{\bf\punct Checker triangulated surfaces and disjoint cycles.}

\begin{observation}
\label{obs:cycles}
For an element of $\wh \kop_N$:

\sm

--- collections of labels at nearIR-farUV vertices in the counterclockwise
cyclic order are disjoint cycles of the permutation $u$;

\sm

---  collections of labels at midIR-farUV vertices in clockwise cyclic order 
are disjoint cycles of $v$;

\sm

---  collections of labels of plus-triangles at nearIR-midUR vertices in counterclockwise cyclic order  are disjoint cycles of 
$vu$; collections of labels of minus-triangles nearIR-midUR  vertices
are disjoint cycles of $uv $.
\end{observation}

This is clear from Fig. \ref{fig:4}.
\hfill $\square$

\sm

{\bf\punct Coloring chambers.} Now let us paint chambers
of surfaces $\cH\in  \kop_N$ in colors $\gimel_j$
in such a way that $m_j$ chambers be painted in color $\gimel_j$.
Denote by $\kop_{\{m_j\}}$ the set of
surfaces obtained in this way.
Notice that we get a necklace in each nearIR-farUV
 vertex (counterclockwise),
and therefore a color surface determines  a neck mess, say 'nearIR-farUV neck mess'. In the same way, we define 'midIR-farUV neck mess', in this case,
we go around vertices clockwise. We also define 'plus nearIR-midIR neck mess',
in this case, we consider counterclockwise necklaces
of plus-triangles around nearIR-midIR vertices. 

\sm

Now let $\lambda$, $\mu$, $\nu$ be three  neck messes having
$m_j$ beads of each color $\gimel_j$.
Denote by
$$
\kop_{\{m_j\}}[\lambda,\mu,\nu]
$$ 
the set of all colored checker triangulated surfaces
such that the corresponding nearIR-farUV, midIR-farUV, 
and plus-nearIR-midIR neck messes are respectively $\lambda$, $\mu$, 
$\nu$. 

Denote by 
$$
\wh\kop_{\{m_j\}}[\lambda,\mu,\nu]
$$ 
the set of surfaces $\in \kop_{\{m_j\}}[\lambda,\mu,\nu]$
whose chambers are labeled by $1$, $2$, \dots, $N$; recall that
 $N=\sum m_j$.  
Clearly, we have
\begin{equation}
\# \wh\kop_{\{m_j\}}[\lambda,\mu,\nu]
=
N! \sum_{\cH^\circ \in \kop_{\{m_j\}}[\lambda,\mu,\nu]}
\frac{1}{\Aut(\cH^\circ)}.
\label{eq:1}
\end{equation}

For a neck mess $\theta$ we denote by $\Gamma[\theta]$ the corresponding
conjugacy class of $S_N$ by the Young subgroup $Y(\{m_j\})$.
Clearly,
\begin{equation}
\label{eq:2}
\#\Gamma[\theta]=\frac{N!}{\#\Aut(\theta)}.
\end{equation}

\begin{observation}
\label{obs:last-correspondence}
 There is a canonical bijection
between $\wh\kop_{\{m_j\}}[\lambda,\mu,\nu]$ and  the set of all
triples 
\begin{equation*}
u\in\Gamma[\lambda],\quad v\in\Gamma[\mu],\quad
w\in \Gamma[\nu]
\end{equation*}
 such that $w=vu$.
 \end{observation}
 
 Indeed, we apply   Observation \ref{obs:cycles} to $u$ and $v$.
\hfill $\square$ 
 
 \sm
 
 {\bf\punct Colored Hurwitz numbers.}
For each $\theta\in S_N/\!/Y(\{m_j\})$ we denote
an element $\gamma(\theta)\in \C\bigl[S_N/\!/Y(\{m_j\}) \bigr]$
by
$$
\gamma(\theta)=\sum_{p\in \Gamma[\theta]} p.
$$

\begin{theorem}
$$
\gamma[\mu]\cdot\gamma[\lambda]=\sum_{\nu\in S_N/\!/Y(\{m_j\})}
c_{\mu\lambda}^\nu\cdot \gamma[\nu],
$$
where the structure constants are given by the following
`colored Hurwitz numbers':
$$
c_{\mu\lambda}^\nu:=(\# \Aut(\nu))
\sum_{\cH^\circ\in \kop_{\{m_j\}}[\lambda,\mu,\nu]}
\frac{1}{\#\Aut \cH^\circ}
.$$
\end{theorem}

 Indeed, by  Observation \ref{obs:last-correspondence},
 $$
c_{\mu\lambda}^\nu= \frac{\# \kop_{\{m_j\}}[\lambda,\mu,\nu]}{\#\Gamma[\nu]}.
 $$ 
We transform this expression applying \eqref{eq:1} and \eqref{eq:2}.
\hfill $\square$

\noindent
\tt
University of Graz,
\\
\vphantom{.}\hfill
Department of Mathematics and Scientific computing;
\\
High School of Modern Mathematics MIPT;
\\
Moscow State University, MechMath. Dept;
\\
University of Vienna, Faculty of Mathematics.
\\
e-mail:yurii.neretin(dog)univie.ac.at
\\
URL: https://www.mat.univie.ac.at/$\sim$neretin/
\\
\phantom{URL:}
https://imsc.uni-graz.at/neretin/index.html

\end{document}